\documentclass{gtart_h}

\def\ifplaintex{\expandafter\ifx\csname documentclass\endcsname\relax}

\def\gtp{{\mathsurround=0pt\it $\cal G\mskip-2mu$eometry \&\ 
$\cal T\!\!$opology $\cal P\!$ublications}}  

\def\recd{{\small Received:\qua\receiveddate\ifx\reviseddate\relax
\else\qquad Revised:\qua\reviseddate\fi\par}} 

\def\lognumber#1{\def\thelognumber{#1}}
\def\volumenumber#1{\def\thevolumenumber{#1}}
\def\volumeyear#1{\def\thevolumeyear{#1}}
\def\papernumber#1{\def\thepapernumber{#1}}
\def\pagenumbers#1#2{\def\startpage{#1}\def\finishpage{#2}}
\def\published#1{\def\publishdate{#1}}

\def\received#1{\def\receiveddate{#1}}

\def\accepted#1{\def\accepteddate{#1}}

\long\def\asciiabstract#1{\long\def\theasciiabstract{#1}}

\let\\\par\let\thelognumber\relax\let\thevolumenumber\relax
\let\thepapernumber\relax\let\thevolumeyear\relax\let\startpage\relax
\let\finishpage\relax\let\publishdate\relax\let\receiveddate\relax
\let\reviseddate\relax\let\accepteddate\relax\let\theasciititle\relax
\let\theasciiauthors\relax
\let\theasciiabstract\relax

\let\theasciiemail\relax

\ifplaintex
\font\logobig=cmssbx10 scaled 3836
\font\logomed=cmssbx10 scaled 2557
\else
\font\logobig=cmssbx10 scaled 4200
\font\logomed=cmssbx10 scaled 2800
\fi

\long\def\makeagttitle{   
\count0=\startpage
\agt\hfill      
\hbox to 45truept{\vbox to 0pt{\vglue -13truept{\logomed A\kern -.37em{\logobig 
T}\kern -.38em G}\vss}\hss}
\break
{\small Volume \thevolumenumber\ (\thevolumeyear)
\startpage--\finishpage\nl
Published: \publishdate}

\vglue .25truein

{\parskip=0pt\leftskip 0pt plus
1fil\def\\{\par\smallskip}{\Large\bf\thetitle}\par\medskip} \vglue
0.05truein

{\parskip=0pt\leftskip 0pt plus 1fil\def\\{\par}{\sc\theauthors}
\par\medskip}%
 
\vglue 0.03truein 

{\small\leftskip 25truept\rightskip 25truept{\bf Abstract}\stdspace\theabstract

{\bf AMS Classification}\stdspace\theprimaryclass
\ifx\thesecondaryclass\relax\else; \thesecondaryclass\fi\par
{\bf Keywords}\stdspace \thekeywords\par}\vglue 7truept

}   

\ifplaintex
\font\phead=cmsl9 scaled 950
\font\pnum=cmbx10 scaled 913
\font\pfoot=cmsl9 scaled 950
\headline{\vbox to 0pt{\vskip -4.5mm\line{\small\phead\ifnum
\count0=\startpage ISSN 1472-2739 (on-line) 1472-2747 (printed)
\hfill {\pnum\folio}\else\ifodd\count0\def\\{ }%
\ifx\theshorttitle\relax\thetitle\else\theshorttitle\fi\hfill{\pnum\folio}
\else\def\\{ and }{\pnum\folio}\hfill\ifx\theshortauthors\relax\theauthors
\else\theshortauthors\fi\fi\fi}\vss}}
\footline{\vbox to 0pt{\vglue 0mm\line{\small\pfoot\ifnum\count0=\startpage
\copyright\ \gtp\hfill\else
\agt, Volume \thevolumenumber\ (\thevolumeyear)\hfill\fi}\vss}}
\else
\font=cmsl9 scaled 1050
\font\lnum=cmbx10 
\font=cmsl9 scaled 1050
\makeatletter
\def\@oddhead{{\small\lhead\ifnum\count0=\startpage ISSN 1472-2739 
(on-line) 1472-2747 (printed)\hfill {\lnum\number\count0}\else\ifodd\count0
\def\\{ }\ifx\theshorttitle\relax \thetitle \else\theshorttitle\fi\hfill
{\lnum\number\count0}\else\def\\{ and }{\lnum\number\count0}
\hfill\ifx\theshortauthors\relax 
\theauthors\else\theshortauthors\fi\fi\fi}}\def\@evenhead{\@oddhead}
\def\@oddfoot{\small\lfoot\ifnum\count0=\startpage\copyright\ \gtp\hfill\else
\agt, Volume \thevolumenumber\ (\thevolumeyear)\hfill\fi}
\def\@evenfoot{\@oddfoot}
\makeatother
\fi
\let\maketitlepage\makeagttitle
\let\makeshorttitle\maketitlepage
\let\maketitle\maketitlepage

\newwrite\gtoutfile
\long\gdef\makeheadfile{  
{\def\\{, }\def\s{ }
\immediate\openout\gtoutfile head.xxx
\immediate\write\gtoutfile{Proxy-for: \ifx\theasciiauthors\relax
\theauthors\else\theasciiauthors\fi\s<\ifx\theasciiemail\relax\theemail\else\theasciiemail\fi>}
\immediate\write\gtoutfile{\noexpand\\}
\immediate\write\gtoutfile{Authors: \ifx\theasciiauthors\relax
\theauthors\else\theasciiauthors\fi}
{\def\\{ }\immediate\write\gtoutfile{Title: \ifx\theasciititle\relax
\thetitle\else\theasciititle\fi}}
\immediate\write\gtoutfile{Subj-class: GT or SG, GR etc}
\immediate\write\gtoutfile{MSC-class: \theprimaryclass\ifx\thesecondaryclass\relax\else, \thesecondaryclass\fi}
\immediate\write\gtoutfile{Journal-ref: Algebr. Geom. Topol. \thevolumenumber\s
(\thevolumeyear) \startpage-\finishpage}
\immediate\write\gtoutfile{Comments: Published by Algebraic and
Geometric Topology at}
\immediate\write\gtoutfile{\s\s\s  http://www.maths.warwick.ac.uk/agt/AGTVol\thevolumenumber/agt-\thevolumenumber-\thepapernumber.abs.html}
\immediate\write\gtoutfile{\noexpand\\}
\immediate\write\gtoutfile{}
\ifx\theasciiabstract\relax
\immediate\write\gtoutfile{\theabstract}\else
\immediate\write\gtoutfile{\theasciiabstract}\fi
\immediate\write\gtoutfile{}
\immediate\write\gtoutfile{\noexpand\\}
\immediate\write\gtoutfile{}
\immediate\closeout\gtoutfile}}  

\def\maketitlepage{\makeagttitle\makeheadfile}
\let\makeshorttitle\maketitlepage
\let\maketitle\maketitlepage

\lognumber{24}
\volumenumber{5}
\volumeyear{2005}
\papernumber{24}
\pagenumbers{563}{576}
\received{1 February 2005} 
\accepted{31 May 2005}
\published{29 June 2005}

\usepackage{amsmath,amssymb}
\usepackage{graphicx,psfrag,verbatim}
\def\psfraga <#1,#2> #3#4{%
\psfrag {#3}{\smash{\rlap{\kern #1 \raise #2\hbox{#4}}}}}

\theoremstyle{plain}
\newtheorem{thm}{Theorem}[section]
\newtheorem{lem}[thm]{Lemma}
\theoremstyle{definition}

\newcommand{\bd}{\partial}
\newcommand{\x}{\times}
\newcommand{\e}{\epsilon}
\newcommand{\R}{\ensuremath{\mathbb{R}}}
\newcommand{\mS}{\ensuremath{\mathcal{S}}}
\newcommand{\mC}{\ensuremath{\mathcal{C}}}

\begin{document}

\title{Knots on a positive template have a
bounded\\number of prime factors}

\author{Michael C. Sullivan}

\address{Department of Mathematics (4408), Southern Illinois 
University\\Carbondale, IL 62901, USA}
\email{msulliva@math.siu.edu}
\urladdr{http://www.math.siu.edu/sullivan}

\begin{abstract}
Templates are branched 2--manifolds with semi-flows used to model 
``chaotic'' 
hyperbolic invariant sets of flows on 3--manifolds. Knotted orbits on a 
template correspond to those in the original flow. Birman and Williams 
conjectured that for any given template the number of prime factors of 
the knots realized would be bounded.
We prove a special case when the template is {\em positive}; 
the general case is now known to be false.
\end{abstract}

\asciiabstract{%
Templates are branched 2-manifolds with semi-flows used to model
`chaotic' hyperbolic invariant sets of flows on 3-manifolds. Knotted
orbits on a template correspond to those in the original flow. Birman
and Williams conjectured that for any given template the number of
prime factors of the knots realized would be bounded.  We prove a
special case when the template is positive; the general case is
now known to be false.}

\primaryclass{37D45} 
\secondaryclass{57M25}
\keywords{Hyperbolic flows, templates, prime knots, composite knots, 
positive braids}
\maketitle

\section{Introduction}

Templates are compact branched 2--manifolds with semi-flows used to model 
certain hyperbolic flows on 3--manifolds. Knotted orbits on a template 
correspond to those in the original flow. Birman and Williams conjectured 
that for any given template
the number of prime factors of the knots realized would be bounded; 
see \cite{BW2}. A counter example was first constructed in \cite{S2}, 
but also see \cite{G}. Here we prove that a for the subclass of 
{\it positive templates} the Birman--Williams conjecture is true.
Section 2 gives background on templates; see also \cite{GHS}. 
Section 3 reviews Cromwell's 
Theorem on factoring positive braids \cite{C}; it is our major tool. 
Some terminology for knots and braids is reviewed, but readers new to
knot theory may want to have the text \cite{BZ} on hand.

\section{Templates}

\begin{figure}[ht!]
	\begin{center}  
        \includegraphics[height=1in]{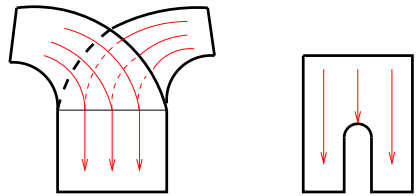}
	\end{center}
	\caption{The charts}
	\label{fig_charts}
\end{figure}

Templates are formed from a finite complex with two types of charts:
{\it joining charts} and {\it splitting charts}, shown in Figure 
\ref{fig_charts}. In the joining charts the flow lines merge at a 
{\it branch line}. There are two entrance segments and one exit segment
in the boundary. The semi-flow is tangent to the rest of the boundary.
The splitting chart has one entrance segment, but its exit set is partitioned
into three sub-segments, indicated by an inward curving of the middle 
sub-segment. The semi-flow is tangent to the two side segments. A template 
is formed by attaching exit sets to entrance sets. It is 
required that in a template the exit set consists of the middle 
portions of the splitting charts and that the entrance set be
empty. It follows that the number of joining charts is equal to the 
number of splitting charts. 

The {\it invariant set} of a template is the set of orbits of the 
semi-flow that never exit. The invariant set is the suspension 
(torus mapping) of a one-sided shift of finite type. (Its inverse 
limit is a suspended two-sided shift of finite type.) Thus, the 
invariant set contains infinitely many closed orbits. In a template 
embedded in $\R^3$ (we will always be working with a given embedding)
the closed orbits form knots. These determine infinitely many 
knot types \cite{FW}; in some cases they support all knot types 
\cite{G}. Franks and Williams \cite{FW} have shown that any 
template can be braided. That is any template can be isotoped so 
that all the closed orbits are presented as braids. If $T$ denotes 
a template we also use $T$ to denote the set of knot realized as 
periodic orbits in the semi-flow.

\begin{figure}[ht!]
	\begin{center}  
        \includegraphics[height=1in]{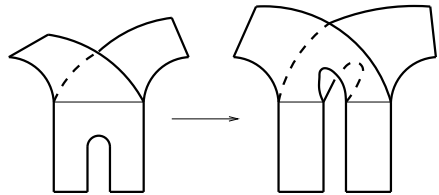}
	\end{center}
	\caption{A split move}
	\label{fig_splitmove}
\end{figure}

We define the {\it split move} via Figure \ref{fig_splitmove}. 
It changes the  topology of a template but does not effect the 
invariant set. 

Knots can be uniquely factored (up to order) in to primes \cite{BZ} 
-- more on this in Section 3. Williams showed that the Lorenz 
template, which has two charts and no twisting in its bands,
contains only prime knots \cite{W}. This and other considerations 
led Birman and Williams \cite{BW2} to conjecture that for any 
given template there would be a  finite bound on the number of 
prime factors for the supported knots. 

Although the original Birman--Williams Conjecture is false work 
in \cite{S2,S3} lead to the weakened conjecture that if a
template had a braid presentation in which all crossings were 
of the same type, then there would be a bound on the number of 
prime factors of the periodic orbits. Such templates are 
called {\it positive templates}. This weakened Birman--Williams
Conjecture is Theorem \ref{thm_main}. Our major tool is a very 
powerful theorem on factoring positive braids due to Peter Cromwell
\cite{C}.

\section{Cromwell's Factoring Theorem}

Let $k$ be a knot, an embedded 1--sphere in $\R^3$. The {\it knot-type} of 
$k$ is its equivalence class under ambient isotopy. An {\em unknot} or 
{\em trivial knot} is any knot equivalent to a circle. 
A projection $\pi$ of $k$ into $\R^2$ is {\it regular} if
the self-intersection set of $\pi(k)$ 
consists of a finite number of transverse double points. 
We say $\pi(k)$ is {\it irreducible} if it has no cut points. 

A knot $k$ is said to be {\it factored} by a 2--sphere $S$ in $\R^3$ 
if $k \cap S$ is transverse and consists of just two points. The factors
are two knots $k_1$ and $k_2$ formed by taking the union of any simple curve 
on $S$ which has as end points $k \cap S$ and the portions of $k$ inside and 
outside of $S$ respectively. If there exists a factoring 2--sphere such that  
neither factor is the unknot then $k$ is a 
{\it composite knot}, and we write $k = k_1 \# k_2$. If the only factors 
of $k$ are itself and the unknot,
then $k$ is {\it prime}, unless $k$ is the unknot. Schubert established 
that nontrivial knots can be factored uniquely into primes, up to order.
An unknot can only be factored into unknots. See \cite{BZ}.

A smooth knot is in {\em braid form} or is {\em braided} if there is an axis 
with respect to which the theta derivative, in cylindrical coordinates 
about the axis, of some parameterization never changes sign. The regular 
projection onto a plane perpendicular to the axis can then be described
symbolically as follows. Let $n$ be the typical number of intersection 
points of the projection and a radius. We say the braid has $n$ strands.
We number the gaps between strands 1 to $n-1$.
Then the integers  $\{ -(n-1),...,-1,1,....n-1 \}$ are used to specify 
the the order of the crossing. They determine a group under concatenation
called the $B_n$ braid group. Thus $111$ in $B_2$ defines a braid with three
positive crossings. 

Let $B_n$ be the $n$--strand braid group. 
A braid is {\em positive} if all 
its crossings are of the same type. All our braids will be positive, so we 
can denote a braid by a word of positive integers.
Let $b = w_1 \cdots w_p \in B_n$, be positive. Then 
$b$ is {\it decomposable} if there exists positive integers $r<n$ and
$q<p$ such that $w_1,\dots,w_{q-1}$ are less then $r$ and $w_q, \dots,w_p$
are greater than or equal to $r$. E.g., $122112234343344$ is decomposable; 
we have $122112234343344$ = $1221122$ \# $12121122$.

\begin{thm} [Cromwell's Theorem] \label{thm_C}
Let $b$ be a positive braid that is an irreducible 
projection of a knot $k$. Then $k$ is prime if and only if 
$b$ is not decomposable.
\end{thm}

Cromwell's approach is to study the intersection of a would-be
factoring sphere with the knot's Seifert surface. In an unpublished 
note \cite{S4} Cromwell's Theorem is proved using a template like 
construction. 
Ozawa \cite{Makoto} extended Cromwell's Theorem to positive knots (knots with 
positive projections, but not necessarily representable as positive braids;
the 5--knot is an example). Ozawa's proof uses incompressible tori and 
is far more elegant than Cromwell's original proof or that given in \cite{S4}.
 
A {\it factoring sphere system} for a composite knot is a 
disjoint set of 2--spheres that factor the knot into primes.
Let $k$ be a positive braid with $n$ prime factors. 
Then it follows from Cromwell's Theorem that there is a factoring 
sphere system for $k$ which consists of $n-1$ concentric spheres
meeting $\R^2$ in $n-1$ concentric circles about the 
braid axis. It will be convenient to allow for small deformations 
in the circles.

\section{The Theorem}

\begin{thm} \label{thm_main}
For any positive braided template $T$, there exists a positive 
integer $N=N(T)$, such that for every knot $k$ in $T$, the number 
of prime factors of $k$ is less than or equal to $N$.
\end{thm}

\proof For a knot $k$ let
$F(k)$ be the number of prime factors of $k$. 
Let $T$ be a positive template with $k$ a closed orbit.
Let $J$ be the number of joining charts, and $B=2J$ be the 
number of bands. Let
$N = 1+\mbox{dim} H_1(T) + J(1+(2J)!)(2J(1+(2J)!)-1)$. 
We will show that  $F(k) \leq N$. 

\begin{figure}[ht!]
	\begin{center}  
        \includegraphics[height=1.5in]{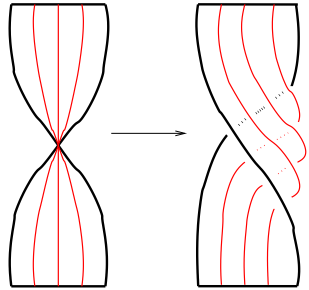}
	\end{center}
	\caption{Projecting a half twist}
	\label{fig_twist}
\end{figure}

We specify a very nice projection of $T$ into the plane. 
Let $\pi:\R^3 \to P \approx \R^2$ be the projection that sends 
$(x,y,z)$ to $(x,y,0)$.
We position $T$ in $\R^3$ so that $T$ is always in $P \x [0,\e]$ 
and $\pi(k)$ has only transverse crossings for any $k \in T$.
(In this paragraph $k$ stands for any periodic orbit of $T$.)
Place each branch line parallel the $x$--axis with the semi-flow 
coming down (decreasing $y$). The bands remain in $P$ with three 
exceptions. 
(i) Where a band has a half twist 
it will go above $P$ ($z>0$) but stay within $P \x [0,\e]$ and the 
$T$ is isotoped so that $\pi(k)$ is transverse, as in 
Figure~\ref{fig_twist}. (ii) When two bands cross we insure 
that $\pi(k)$ is transverse,
$T \subset P \x [0,\e]$, and we do not allow more than two bands to 
cross at a time. (iii) Just above (in the $y$ direction) each branch 
line we insure $\pi(k)$ is transverse and $T \subset P \x [0,\e]$; 
see again the joining chart in Figure \ref{fig_charts}.

Suppose $k \in T$ is composite. By Cromwell's Theorem there exists a 
collection of concentric topological circles in $P$,
$\mC = \{ C_1,...,C_n \}$ that factor $\pi(k)$ into  $n+1=F(k)$ 
primes; assume $C_i$ is interior to $C_{i+1}$ for $i=1,...,n-1$. 
Let $\mS = \{ S_i = C_i \x [-i, i] 
                       \:\mbox{union two disks}\: \}_{i=1}^{n}$.
Then $\mS$ is a complete factoring sphere system for $k$ as 
it appears in $T$. We isotope the $C_i$'s so that $\mS$ is 
transverse to $T$ and is still a complete factoring sphere system. 
Hence, $\mS \cap T$ is compact.  

Let $\beta$ be the set of branch points of $T$. The intersection 
$\mS \cap T$ determines a finite 1--complex where the vertices 
are the points $\mS \cap (\bd T \cup \beta)$. The points 
$\mS \cap \bd T$ have valence one, while the points 
$\mS \cap \beta$ have valence three.

The one-dimensional simplices of $\mS \cap T$ are classified as 
follows. The boundary of $T$ can be partitioned into segments 
(one-simplices) whose end points are on the branch lines. Denote 
by $\bd^0 T$ the union of those segments of $\bd T$ where the 
semi-flow never exits. Call the remaining segments 
{\em band splitting} or {\em exit} segments.

\begin{itemize}

\item $\lfloor$--segments connect a branch point to a point on 
      $\bd^0 T$ below and to the right.

\item $\rfloor$--segments connect a branch point to a point on 
      $\bd^0 T$ below and to the left.

\item $\lceil$--segments connect a branch point to a point on $\bd T$ 
      above and to the right.

\item $\rceil$--segments connect a branch point to a point on $\bd T$ 
      above and to the left.

\item $bs$--segments connect a branch point to a point on an exit 
      segment.
 
\item $bb$--segments, or {\em branch-to-branch} segments connect one 
      branch line to another.

\item $ss$--segments, or {\em edge-to-edge} segments connect one 
      side of a band to the other. 

\item $\cap$--segments connect two points on a branch line from above.

\item $\cup$--segments connect two points on a branch line from below. 

\item $($--segments connect two points of a segment of $\bd T$.

\item Trivial loops are loops that miss $\bd T \cup \beta$.

\end{itemize}

The choice of $\mC$ and hence $\mS$ is far from unique. We shall 
insist on the following minimality assumptions.

\begin{itemize}

\item The number of segments in $\mS \cap T$ is the smallest among
all prime factoring sphere systems of $k$, as constructed above.

\item The number of branch points in $\mS \cap T$ is the smallest 
possible relative to the assumptions above.
\end{itemize}

\begin{lem}
There are no trivial loops, $($--, $\lfloor$--, $\rfloor$--, $\cup$-- or 
$\cap$--segments in $\mS \cap T$.
\end{lem}

\begin{figure}[ht!]\small
        \psfrag{a}{a}\psfrag{b}{b}\psfrag{c}{c}
	\begin{center}  
        \includegraphics[width=5in]{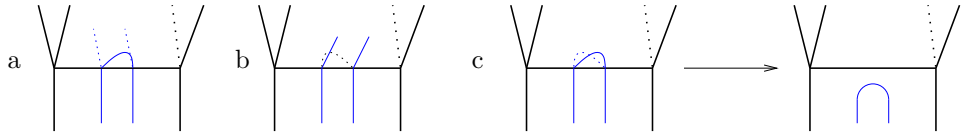}
	\end{center}
	\caption{No $\cap$--segments}
	\label{fig:nocaps}
\end{figure}

\begin{proof}
If a trivial loop in $\mS \cap T$ meets the knot $k$ then a trivial 
factor is produced. This is not permitted. If a trivial loop misses $k$ 
the corresponding sphere misses $k$. This is not permitted. 
If a $($--segment meets $k$ a trivial factor in produced. If a 
$($--segment misses $k$ we may assume it is inner most and 
deform the corresponding sphere to eliminate it, reducing the 
number a segments in $\mS \cap T$.

For $\lfloor$--, $\rfloor$--, and $\cup$--segments the arguments are similar 
and can be found in Lemma 1.1 of \cite{W}.

For $\cap$--segments we consider three cases, (a), (b) and (c) as shown in 
Figure~\ref{fig:nocaps}. In (a) and (b) $\pi(\mS \cap T)$ has valence 
three 
points contradicting the fact that $\mC$ consists of a union of simple 
closed curves. (We shall say that there are no Y's in $\pi(\mS \cap T)$.)
The configuration in (c) can be deformed to yield a factoring
sphere system $\mS'$ with fewer segments in $\mS' \cap T$, contradicting 
the minimality assumptions. 
\end{proof}

\begin{lem}
The connected components of $\mS \cap T$ consist of three types.
\begin{enumerate}
\item $ss$--segments.
\item Nontrivial trees (trees with more than one segment).
\item Graphs consisting of a single cycle and some $\lceil$-- or 
      $\rceil$--segments.
\end{enumerate}
\end{lem}

\begin{proof}
If a component contains an $ss$--segment, it is an $ss$--segment.
We need show that a non-tree component has only one cycle. 
A cycle is formed only from $bb$--segments. Pick a vertex on a cycle. 
Tracing down (with the semi-flow) we must encounter a $bb$--segment. 
If there are two $bb$--segments above our vertex, then there is a Y is 
the projection. Thus, above the vertex there is one $bb$--segment and 
either a $\lceil$-- or $\rceil$--segment.
\end{proof}

\begin{figure}[ht!]\small
	\begin{center}  
        \psfrag{Y}{Y}
        \psfrag{push}{push}
        \psfrag{down}{down}
        \includegraphics[height=1.8in]{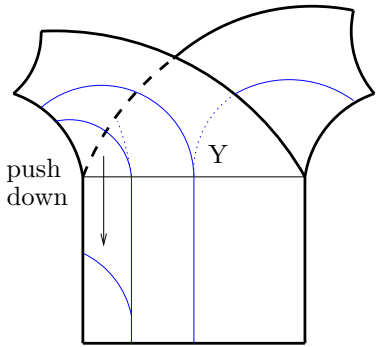}
	\end{center}
	\caption{Treetop analysis}
	\label{fig_treetop}
\end{figure}

Let $\tau$ be a non-trivial tree component in $\mS \cap T$. Pick a 
point on $\tau$ and trace down (with the flow direction). This path 
must exit the template somewhere. Since there are no $\rfloor$-- or 
$\lfloor$--segments there must be a $bs$--segment. (The trunk of a 
tree is rooted at a split.) Now trace up. When we meet a branch 
line we make a choice as to which segment to take. If possible we 
avoid $\rceil$-- and $\lceil$--segments in favor of a $bb$--segment. 
But, this path too must terminate. Therefore there is a branch 
line meeting $\tau$ where both of the segments above are $\rceil$-- 
or $\lceil$--segments. Call this the {\em treetop}. If both segments 
are the same type, minimality is violated; pushing the sphere down 
through the branch line reduces the number of segments. 
Furthermore, the front one must be a $\lceil$--segment, and the back
one must be a $\rceil$--segment or else the projection $\pi(\tau)$ 
will contain a Y. See Figure \ref{fig_treetop}.

\begin{figure}[ht!]
	\begin{center}  
        \includegraphics[height=3in]{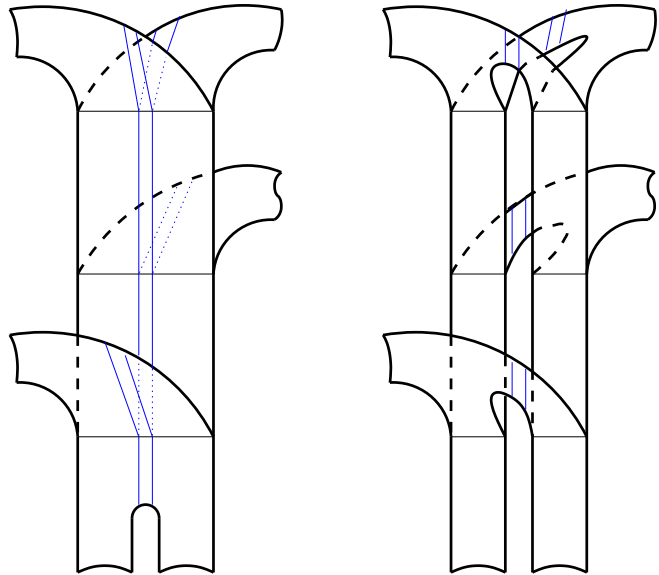}
	\end{center}
	\caption{Removing a pair of parallel trees}
	\label{fig_treecut}
\end{figure}

We use split moves to remove all tree components from $\mS \cap T$. 
Figure~\ref{fig_treecut} gives an example. Of course each split move 
changes the number charts and bands; $J$ will increase by 1, and 
$B$ will increase by 2. (The number of $ss$--segments is also 
increased, by 2 at a tree-top, and by 1 otherwise, per tree.) An 
upper bound on the number of split moves needed to remove all of 
the tree components for a given template $T$ can be derived from 
the following facts.
\begin{enumerate}
\item Since a tree component projects into a braided circle it 
      cannot meet the same branch line twice.
\item ``Parallel'' trees, those that use the same bands, are removed 
      by the same sequence of split moves. See 
      Figure~\ref{fig_treecut}.
\item The number of sets of parallel groupings of trees is bounded 
      by $B!$. This follows from (1).
\end{enumerate}
Therefore, the maximum number of split moves needed to remove all the
tree components from $T$ is $JB!$. Call the new template formed $T'$.

Notation: For an oriented knot $k$ and points $a$ and $b$ on $k$,
let $k(a,b)$ denote the oriented arc in $k$ starting at $a$ and 
ending at $b$.

\begin{lem}
The number of non-tree components is bounded by $\dim H_1(T)$.
\end{lem}

\begin{proof}
Any cycle in $\mS \cap T'$ corresponds to a cycle in $\mS \cap T$ 
since the split moves do not introduce new non-tree components. We 
will show that if  $\mS \cap T'$ has two parallel cycles, by which 
we mean they pass through the same bands (they are homologous), then 
the factoring of $k$ by $\mS$ produces an unknotted factor. This 
contradiction will give the result.

Suppose two graphs components ($G_1$, $G_2$) have parallel cycles 
($C_1$, $C_2$). Assume the cycles are inner most, that is there are 
no other cycles in between them. Thus they bound an annulus $A$ in 
$T'$. Let $S_1$ and $S_2$ be the corresponding spheres with  $S_1$  
inside $S_2$. The annulus $A$ meets only these two spheres and 
$A \cap \mS = C_1 \cup C_2$. The two spheres partition $\R^3$ into 
three regions with the interior of $A$ in between $S_1$ and $S_2$.

The knot $k$ pierces each sphere ($S_1$, $S_2$) exactly twice. 
First suppose the knot $k$ pierces each component ($G_1$, $G_2$) 
exactly twice. Let $k \cap G_i = \{p_i, q_i\}$, for $i=1,2$. We will 
construct a closed loop $u$ that is a factor of $k$. We will show 
that $u$ is an unknot, thus deriving a contradiction. Start from 
$p_2$ and assume without loss of generality that $k$ passes from 
outside $S_2$ to its inside. 

If $k$ meets $q_2$ before hitting $S_1$ form $u$ by uniting the arc 
of $k$ from $p_2$ to $q_2$ with a circular arc in $C_1$. Since $u$ 
is embedded in an annulus it is unknotted.

Assume $k$ enters $S_1$ at $p_1$ and re-emerges at $q_1$. Its next 
intersection with $\mS$ will be at $q_2$. Form $u$ by taking the 
union of the arc $k(p_2,p_1)$, an arc of $G_1$ connecting $p_1$ 
to $q_1$, the arc $k(q_1,q_2)$, and an arc of $G_2$ connecting 
$q_2$ to $p_2$.  We chose the arcs in $G_1$ and $G_2$ so that $u$ 
is braided (although we may need to make a small isotopy if these 
arcs start or end on $\lceil$-- or $\rceil$--segments.) We divide 
the problem into subcases.

First suppose $k$ misses the annulus $A$. Thus, $p_1$, $p_2$, $q_1$, 
and $q_2$ are in $\lceil$-- or $\rceil$--segments. The 
construction for $u$ is shown in Figure~\ref{fig_case1} where $u$ 
is seen to be an unknot.

\begin{figure}[ht!]\small
	\begin{center} 
	\psfrag{s1}{$S_1$}
	\psfrag{s2}{$S_2$}
	\psfraga <-2pt, 0pt> {p1}{$p_1$}
	\psfrag{p2}{$p_2$}
	\psfrag{q1}{$q_1$}
	\psfrag{q2}{$q_2$}
	\psfrag{k}{$k$}
	\psfrag{u}{$u$}
        \includegraphics[height=2.5in]{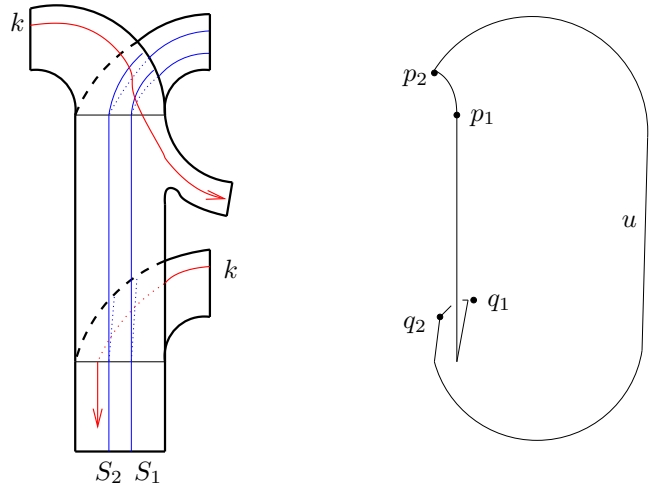}
	\end{center}
	\caption{$u$ is unknotted.}
	\label{fig_case1}
\end{figure}

Now suppose $p_2$ is in a $\rceil$--segment, but that $k$  enters 
$A$, say at a point $x$ on the branch line $\beta_1$. Our $k$ may 
wind around on $A$ but will not meet $\beta_1$ to the left of $x$, 
otherwise it could not get to $G_1$. Suppose it winds around $m-1$ 
times before meeting $G_1$. (Notice $p_1 \in C_1$.) So far $u$ has 
no crossings. When $k$ re-emerges from $S_1$ it must do so through a 
$\lceil$--segment, otherwise it cannot get back to $G_2$. If $k$ exits 
$S_2$ without meeting $A$ at a $\lceil$--segment, $u$ will have 
braid word $m(m-1)\cdots 2 1$, which is an unknot. 
If $k$ does re-enter $A$, say at a point $y$ of the branch line 
$\beta_2$, then $y$ is to the left of every point of 
$k(p_2,p_1) \cap \beta_2$, otherwise $k$ will never get back to $G_2$. 
Suppose $k(q_1,q_2)$ wraps around $A$ $n-1$ times before exiting $S_2$.
Then the braid word of $k$ is of the form $12\cdots n m (m-1) \cdots 
(n+1)$; see Figure~\ref{fig_case2}. Again, $u$ is an unknot.

\begin{figure}[ht!]\small
	\begin{center}  
        \psfraga <-1.5pt,0pt> {x}{$x$}
        \psfraga <-1.5pt,0pt> {y}{$y$}
        \psfraga <-1pt,0pt> {b1}{$\beta_1$}
        \psfraga <-3pt,0pt> {b2}{$\beta_2$}

        \includegraphics[height=2.7in]{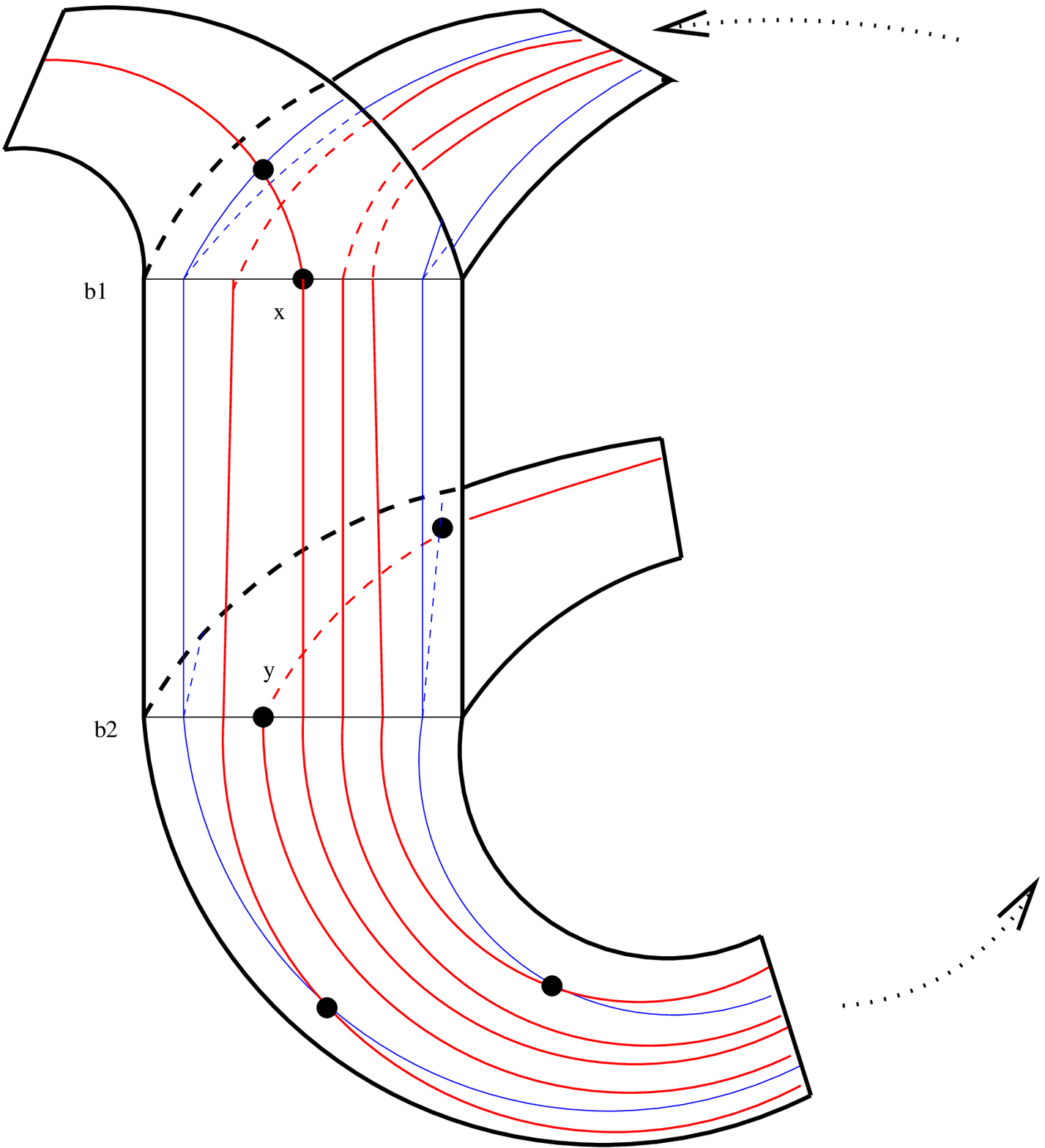}

	\vspace{.3in}

        \includegraphics[height=2.7in]{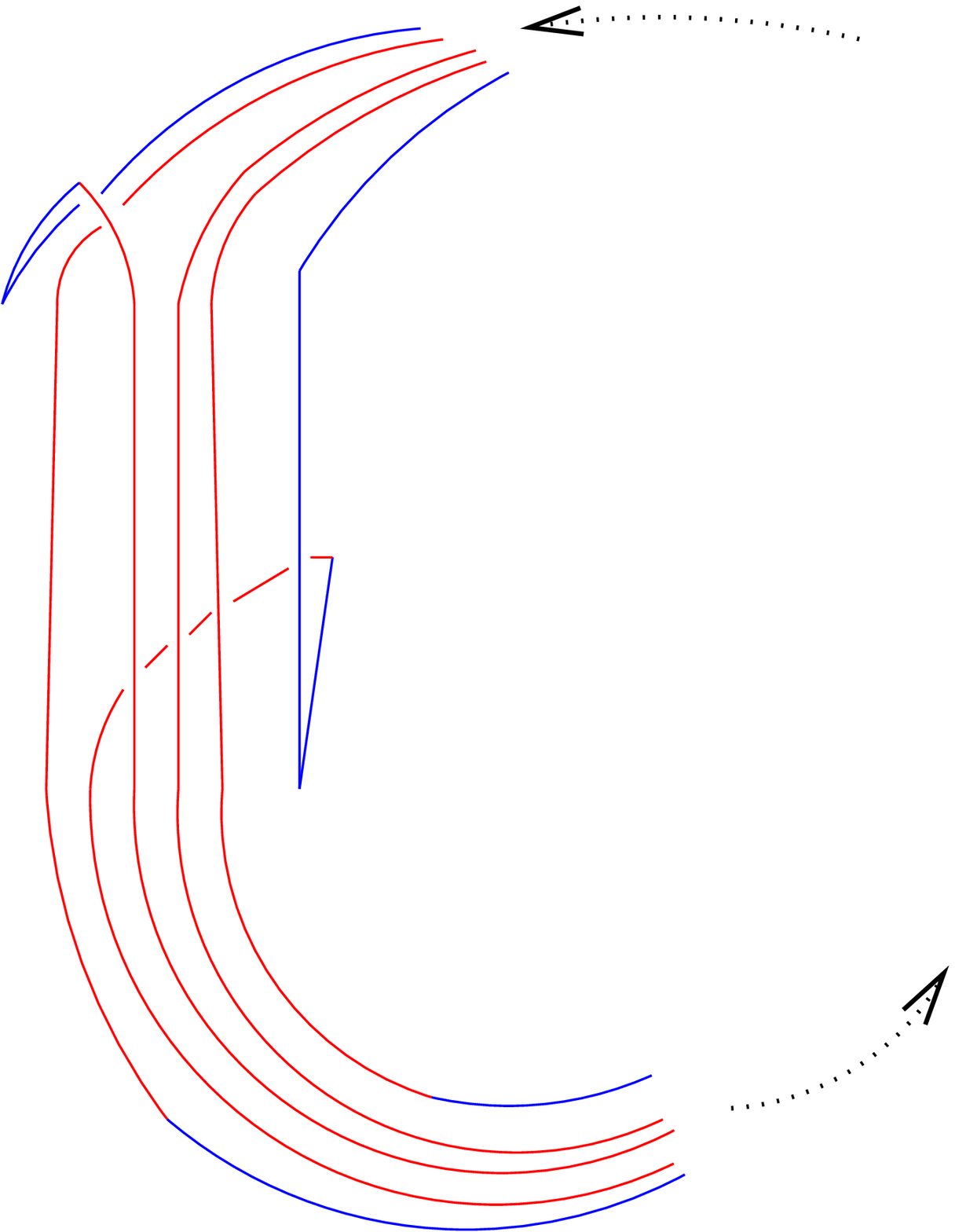}
	\end{center}
	\caption{$u$ is unknotted}
	\label{fig_case2}
\end{figure}

It may be that $k$ misses both $G_1$ and $G_2$. In this case, $p_1$,
$p_2$, $q_1$, and $q_2$ are on $ss$--segments just below and above 
bands that $C_1$ and $C_2$ pass through. The construction of $u$ is 
now very similar to the subcase above where $k$ missed the annulus 
$A$. The only difference is that $u$ will have small segments on 
the spheres that are outside of the template. See Figure 
\ref{fig_case3}.

\begin{figure}[ht!]
	\begin{center}  
        \includegraphics[height=2.5in]{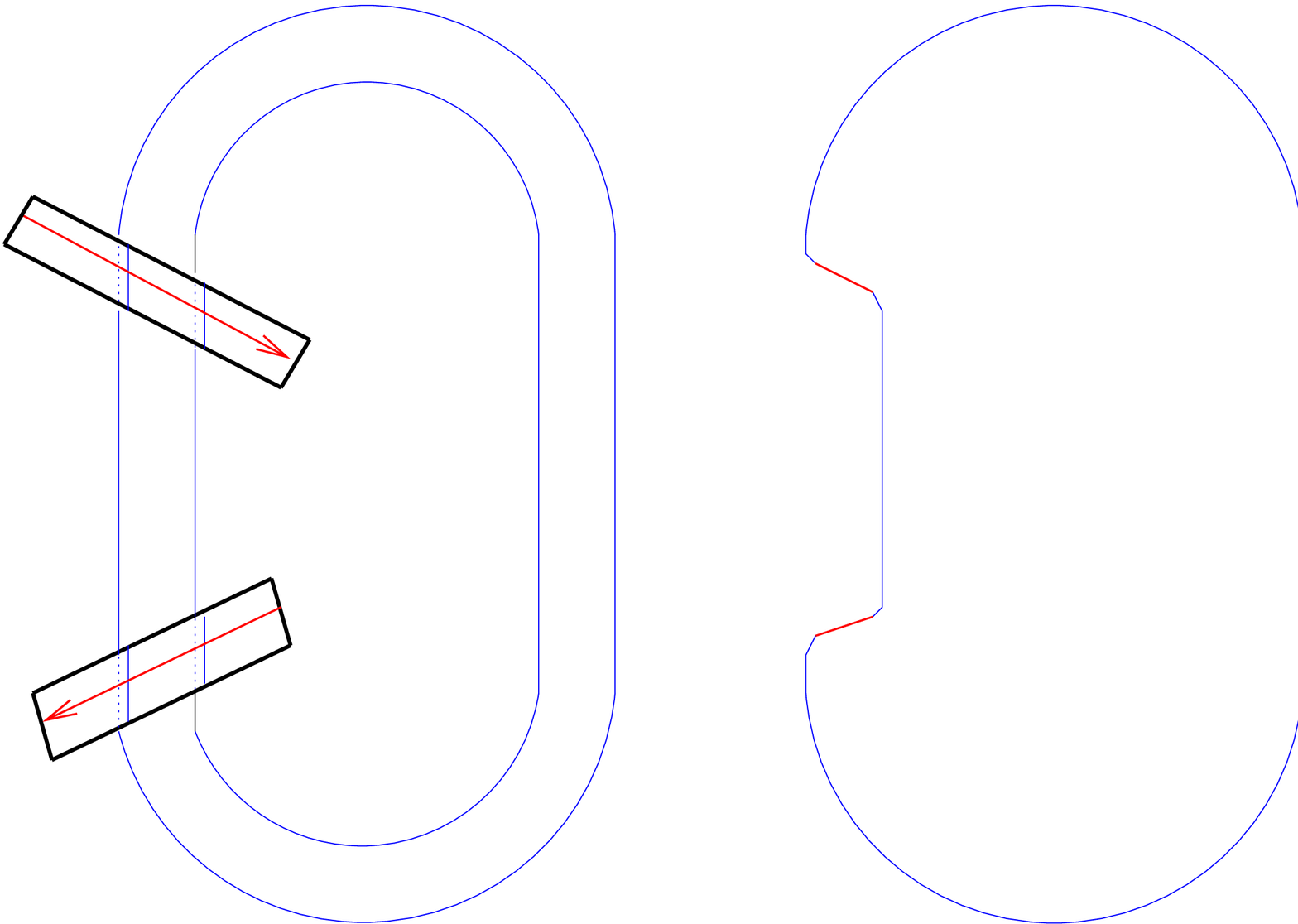}
	\end{center}
	\caption{$u$ is unknotted}
	\label{fig_case3}
\end{figure}

Now if $k$ enters $S_1$ and $S_2$ through $G_1$ and $G_2$ but exits
through 
$ss$--segments (or vise versa) it is not hard to show that $u$ 
will have braid word of the form $m\cdots 1$ and is thus unknotted.
There are no other cases.
\end{proof}

We now turn to the $ss$--segments in $T'$. Our goal is to bound the 
number of spheres needed to factor knots in $T$. Thus we only need 
to bound the number of $ss$--segments that meet $k$. We classify 
such $ss$--segments into two subtypes. Let $E$ be an $ss$--segment 
that meets $k$. Let $C_i$ be the circle containing $E$. At some 
point $p$, $k$ meets $C_i \cap T'$ once again. If the component 
that $p$ is in is a non-tree graph, call $E$ an $ssg$--segment. If 
the component that $p$ is in is another $ss$--segment, call $E$ an 
$ssss$--segment. In this case the $ss$--segment containing $p$ is 
denoted $\hat{E}$, and $E$ and $\hat{E}$ are called {\em associated} 
$ssss$--segments. Clearly, the number of $ssg$--segments is bounded 
by the number of non-tree components. Let $B'$ be the number of 
bands in $T'$; $B' \leq B + 2J(2J)!=B(1+B!)$, since each split 
move produces two additional bands.

\begin{lem}
The number of $ssss$--segments in $T'$ is bounded by $B'(B'-1)$.
\end{lem}

\begin{proof}
The proof is divided into two claims. 

{\bf Claim 1}\qua {\sl Two associated $ssss$--segments cannot be 
in the same band.}
 
Let $E$ and $\hat{E}$ be associated $ssss$--segments. 
Assume they are inner most among such pairs. 
There cannot be a segment from another circle between
them. If the knot misses  $E$ and $\hat{E}$ we can
deform the sphere so as to reduce the number of 
segments by two. If the knots meet  $E$ or $\hat{E}$,
it meets both, and a trivial factor is produced,  
as Figure~\ref{fig:EE} shows.

\begin{figure}[ht!]\small
	\begin{center}  
	\psfrag{E1}{$E$}
	\psfrag{E2}{$\hat{E}$}
	\psfrag{Ci}{$C_i$}
        \includegraphics[height=1in]{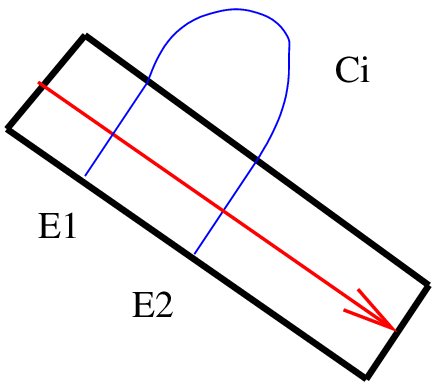}
	\end{center}\vspace{-5mm}
	\caption{A trivial factor}
	\label{fig:EE}
\end{figure}

Let $E_1$ and $E_2$ be $ssss$--segments in the same band $b$.
Let $\hat{E}_1$ and $\hat{E}_2$ be their respective associates.

{\bf Claim 2}\qua {\sl The segments $\hat{E}_1$ and $\hat{E}_2$ cannot be 
in the same band.}

Suppose they were both in the band $\hat{b}$. We can assume such a 
pairing is inner most. Then Figure~\ref{fig:bb} shows that a trivial 
factor would be produced.

\begin{figure}[ht!]\small
	\begin{center}  
	\psfrag{E1}{$E_1$}
	\psfrag{E2}{$E_2$}
	\psfrag{E1'}{$\hat{E}_1$}
	\psfrag{E2'}{$\hat{E}_2$}
        \includegraphics[height=1in]{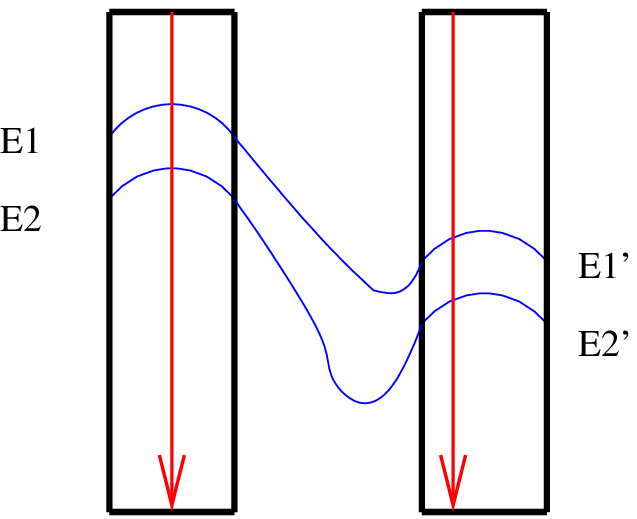}
	\end{center}\vspace{-5mm}
	\caption{A trivial factor}
	\label{fig:bb}
\end{figure}
 
The two Claims imply the desired bound holds.
\end{proof}

The lemmas above establish that $F(k) \leq N(T)$.


\Addresses\recd

\begin{thebibliography}{100}

\bibitem{BW2}
J. Birman \& R. Williams. Knotted Periodic Orbits in Dynamical 
Systems II: Knot Holders for Fibered Knots. {\em Contemporary Math.} 
{\bf 20} (1983) 1--60.
  \MR{0718132}

\bibitem{BZ}
G. Burde \& H. Zieschang. Knots,
Second edition, de Gruyter Studies in Mathematics, 5,
Walter de Gruyter \& Co., Berlin, 2003.
  \MR{1959408}

\bibitem{C}
P. Cromwell. 
Positive braids are visually prime.
{\em Proc. London Math. Soc.} (3) {\bf 67} (1993) 384--424.
  \MR{1226607}

\bibitem{FW}
J. Franks \& R. Williams.
Entropy and knots. 
{\it Trans. Amer. Math. Soc.}  291  (1985),  no. 1, 241--253.
  \MR{0797057}

\bibitem{G}
R. Ghrist.
Branched two-manifolds supporting all links. 
{\em Topology}  36  (1997),  no. 2, 423--448.
  \MR{1415597}

\bibitem{GHS}
R. Ghrist, P. Holmes \& M. Sullivan.
{\em Knots and links in Three-Dimensional Flows},
Lecture Notes in Mathematics, Vol. 1654,
Springer-Verlag, Berlin, 1997.
  \MR{1480169}

\bibitem{Makoto}
M. Ozawa.
Closed incompressible surfaces in the complements of positive knots. 
{\em Comment. Math. Helv.} 77 (2002), no. 2, 235--243.
  \MR{1915040}

\bibitem{S2}
M. Sullivan. 
Composite knots in the Figure-8 knot complement can have any number 
of prime factors, 
{\em Top. and its Appl.} 55 (1994) 261--272.
  \MR{1259509}

\bibitem{S3}
M. Sullivan.
The prime decomposition of knotted periodic orbits in dynamical 
systems,
{\em The Journal of Knot Theory and its Ramifications}, 
Vol. 3 No. 1 (1994) 83-120.
  \MR{1265454}

\bibitem{S4}
M. Sullivan.
Factoring positive braids via branched manifolds.
Preprint. 
\url{http://galileo.math.siu.edu/~msulliva/Preprints/}

\bibitem{W}
R. F. Williams. Lorenz Knots are Prime,
{\em Ergod. Th. \& Dynam. Sys.} {\bf 4} (1983) 147--163.
  \MR{0758900}

\end{thebibliography}
\end{document}